\documentclass[12pt]{article}
\usepackage{amsfonts}
\usepackage{latexsym,amsmath,amssymb}
\usepackage{graphics,epstopdf}
\usepackage[pdftex]{graphicx}
\usepackage{subfig}

\numberwithin{equation}{section}
\begin{document}

\bigskip

\bigskip

\begin{center}
{\Large \textbf{\ Some approximation results for Bernstein-Kantorovich
operators based on $(p,q)$-calculus }}

\bigskip

\bigskip

\textbf{M. Mursaleen}, \textbf{Khursheed J. Ansari} and \textbf{Asif Khan}

Department of\ Mathematics, Aligarh Muslim University, Aligarh--202002, India%
\\[0pt]

mursaleenm@gmail.com; ansari.jkhursheed@gmail.com; asifjnu07@gmail.com \\[0pt%
]

\bigskip

\bigskip

\textbf{Abstract}
\end{center}

\parindent=8mm {\footnotesize {In this paper, we introduce a new analogue of
Bernstein-Kantorovich operators which we call as $(p,q)$%
-Bernstein-Kantorovich operators. We discuss approximation properties for
these operators based on Korovkin's type approximation theorem and we
compute the order of convergence using usual modulus of continuity and also
the rate of convergence when the function $f$ belongs to the class Lip$%
_M(\alpha)$. Moreover, we also study the local approximation property of the
sequence of positive linear operators $K_n^{(p,q)}$. We show comparisons and
some illustrative graphics for the convergence of operators to a function.
In comparison to $q$-analogoue of Bernstein-Kantorovich operators, our
generalization gives more flexibility for the convergence of operators to
a function. %Furthermore, we give comparisons and some
%illustrative graphics for the convergence of operators to some function.
}}

\bigskip



\section{Introduction and preliminaries}

\parindent8mm During the last two decades, the applications of $q$-calculus
emerged as a new area in the field of approximation theory. The rapid
development of $q$-calculus has led to the discovery of various
generalizations of Bernstein polynomials involving $q$-integers. The aim of
these generalizations is to provide appropriate and powerful tools to
application areas such as numerical analysis, computer-aided geometric
design and solutions of differential equations.

\parindent8mmUsing $q$-integers, Lupa\c{s} \cite{lp} introduced the first $q$%
-Bernstein operators \cite{brn} and investigated its approximating and
shape-preserving properties. Another $q$-analogue of the Bernstein
polynomials is due to Phillips \cite{pl}. Since then several generalizations
of well-known positive linear operators based on $q$-integers have been
introduced and studied their approximation properties. For instance, $q$%
-Bleimann, Butzer and Hahn operators \cite{ar1}; $q$-parametric Sz\'{a}%
sz-Mirakjan operators \cite{m1}; $q$-Bernstein-Durrmeyer operators \cite{m2}%
; $q$-analogue of Sz\'{a}sz-Kantorovich operators \cite{m3}.

Recently, Mursaleen \textit{et al} introduced $(p,q)$-calculus in
approximation theory and constructed the{\ $(p,q)$-analogue of Bernstein
operators \cite{mka1} and $(p,q)$-analogue of Bernstein-Stancu operators
\cite{mka2}, $(p,q)$-analogue of Bleimann-Butzer-Hahm operators \cite{mnak},
Bernstein-Schurer operarors \cite{mnn} and investigated their approximation
properties. }The $(p,q)$-analog of Sz\'{a}sz-Mirakyan operators \cite{acar},
Kantorovich type Bernstein-Stancu-Schurer operators \cite{cai} and
Kantorovich variant of $(p,q)$-Sz\'{a}sz-Mirakjan operators \cite{maa} have
recently been studied too.

\parindent=8mm Motivated by their work, in this article, authors introduce a
new analogue of Bernstein-Kantorovich operators. Paper is organized as
follows: In Section 2, we define $(p,q)$-Bernstein-Kantorovich operators and
establish a bsic lemma which is used in proving main results. In Section 3,
we discuss approximation properties for these operators based on Korovkin's
type approximation theorem and we compute the order of convergence using
usual modulus of continuity and also the rate of convergence when the
function $f$ belongs to the class Lip$_{M}(\alpha )$. Moreover, we also
study the local approximation property of the sequence of positive linear
operators $K_{n}^{(p,q)}$. In Section 4, we give some examples to show
comparisons and some illustrative graphics for the convergence of operators
to a function.

Let us recall certain definitions and notations of $(p,q)$-calculus:

The $(p,q)$-integer was introduced in order to generalize or unify several
forms of $q$-oscillator algebras well known in the earlier physics
literature related to the representation theory of single parameter quantum
algebras \cite{chak}. The $(p,q)$-integer $[n]_{p,q}$ is defined by

\begin{equation*}
\lbrack n]_{p,q}:=\frac{p^{n}-q^{n}}{p-q},~~~n=0,1,2,\cdots ,~~0<q<p\leq 1.
\end{equation*}%
\begin{equation*}
\left[ k\right] _{p,q}!:=\left\{
\begin{array}{lll}
\left[ k\right] _{p,q}\left[ k-1\right] _{p,q}...1 & , & k\geq 1, \\
1, &  & k=0%
\end{array}%
\right.
\end{equation*}

The $(p,q)$-Binomial expansion is
\begin{equation*}
(x+y)_{p,q}^n:=(x+y)(px+qy)(p^2x+q^2y)\cdots(p^{n-1}x+q^{n-1}y)
\end{equation*}

and the $(p,q)$-binomial coefficients are defined by

\begin{equation*}
\left[
\begin{array}{c}
n \\
k%
\end{array}%
\right] _{p,q}:=\frac{[n]_{p,q}!}{[k]_{p,q}![n-k]_{p,q}!}.
\end{equation*}

The definite integrals of the function $f$ are defined by
\begin{equation*}
\int_{0}^{a}f(x)d_{p,q}x=(q-p)a\sum\limits_{k=0}^{\infty}\frac{p^k}{q^{k+1}}%
f\left(\frac{p^k}{q^{k+1}}a\right),~~~\text{when}~~~\left|\frac pq\right|<1,
\end{equation*}
and
\begin{equation*}
\int_{0}^{a}f(x)d_{p,q}x=(p-q)a\sum\limits_{k=0}^{\infty}\frac{q^k}{p^{k+1}}%
f\left(\frac{q^k}{p^{k+1}}a\right),~~~\text{when}~~~\left|\frac pq\right|>1.
\end{equation*}

Details on $(p,q)$-calculus can be found in \cite{mah,jacob,jag,vivek,sad}.
For $p=1$, all the notions of $(p,q)$-calculus are reduced to $q$-calculus
\cite{kac}.

\section{Construction of Operators}

Mursaleen et. al \cite{mka1} introduced $(p,q)$-analogue of Bernstein operators as
\begin{equation*}
B_{n,p,q}(f;x)=\sum\limits_{k=0}^{n}\left[
\begin{array}{c}
n \\
k%
\end{array}%
\right] _{p,q}x^{k}\prod\limits_{s=0}^{n-k-1}(p^{s}-q^{s}x)~~f\left( \frac{%
[k]_{p,q}}{[n]_{p,q}}\right) ,~~x\in \lbrack 0,1].
\end{equation*}
But $B_{n,p,q}(1;x)\ne1$ for all $x\in[0,1]$. Hence, they re-introduced their operators in \cite{mka11} as follows :
\begin{equation*}
B_{n,p,q}(f;x)=\frac1{p^{\frac{n(n-1)}2}}\sum\limits_{k=0}^{n}\left[
\begin{array}{c}
n \\
k%
\end{array}%
\right] _{p,q}p^{\frac{k(k-1)}2}x^{k}\prod\limits_{s=0}^{n-k-1}(p^{s}-q^{s}x)~~f\left( \frac{%
[k]_{p,q}}{p^{k-n}[n]_{p,q}}\right) ,~~x\in \lbrack 0,1].
\end{equation*}

The same problem was occurring with the operators introduced in \cite{mka2}. So the revised form of $(p,q)$-analogue of Bernstein-Stancu operators are given as follows:

\begin{equation*}
S_{n,p,q}(f;x)=\frac1{p^{\frac{n(n-1)}2}}\sum\limits_{k=0}^{n}\left[
\begin{array}{c}
n \\
k%
\end{array}%
\right] _{p,q}p^{\frac{k(k-1)}2}x^{k}\prod\limits_{s=0}^{n-k-1}(p^{s}-q^{s}x)~~f\left( \frac{%
{p^{n-k}[k]_{p,q}+\alpha}}{[n]_{p,q}+\beta}\right) ,~~x\in \lbrack 0,1].
\end{equation*}

Note that for $p=1$, $(p,q)$-Bernstein operators and  $(p,q)$-Bernstein-Stancu operators turn out to
be $q$-Bernstein operators and $q$-Bernstein-Stancu operators, respectively.

Dalmanoglu \cite{dal} defined the Bernstein-Kantorovich \cite{kan} operators
using $q$-calculus as follows:

\begin{equation*}
K_{n,q}(f;x)=[n+1]_{q}\sum\limits_{k=0}^{n}p_{n,k}(q;x)\int_{{[k]_{q}}/{%
[n+1]_{q}}}^{{[k+1]_{q}}/{[n+1]_{q}}}f(t)d_{q}t ,~~x\in \lbrack 0,1],
\end{equation*}
\begin{equation*}
p_{n,k}(q;x):= \left[
\begin{array}{c}
n \\
k%
\end{array}%
\right] _{q}x^{k}\prod\limits_{s=0}^{n-k-1}(1-q^{s}x).
\end{equation*}
where $K_{n,q}:$ $C[0,1]\rightarrow C[0,1]$ are defined for any $n\in
\mathbb{N}$ and for any function $f\in C[0,1].$

Now, we introduce $(p,q)$-analogue of Bernstein-Kantorovich operators as
\begin{equation*}
K_{n}^{(p,q)}(f;x)=\frac{[n]_{p,q}}{p^{\frac{n(n-1)}{2}}}\sum%
\limits_{k=0}^{n}\frac{b_{n,k}^{(p,q)}(x)}{p^{n-k}q^{k}}\int_{\frac{[k]_{p,q}%
}{p^{k-n-1}~[n]_{p,q}}}^{\frac{[k+1]_{p,q}}{p^{k-n}~[n]_{p,q}}%
}f(t)d_{p,q}t,~~x\in \lbrack 0,1]\eqno(2.1)
\end{equation*}%
where
\begin{equation*}
b_{n,k}^{(p,q)}(x)=\left[
\begin{array}{c}
n \\
k%
\end{array}%
\right] _{p,q}~(x)_{p,q}^{k}(1-x)_{p,q}^{n-k}=\left[
\begin{array}{c}
n \\
k%
\end{array}%
\right] _{p,q}~p^{\frac{k(k-1)}{2}}x^{k}(1-x)_{p,q}^{n-k}=\left[
\begin{array}{c}
n \\
k%
\end{array}%
\right] _{p,q}p^{\frac{k(k-1)}{2}}x^{k}\prod%
\limits_{s=0}^{n-k-1}(p^{s}-q^{s}x)
\end{equation*}%
and $(x)_{p,q}^{k}:=x(px)(p^{2}x)\cdots (p^{k-1}x)=p^{\frac{k(k-1)}{2}}x^{k}$, and $f$ is a non-decreasing function.\newline

For $p=1$, operators (2.1) turns out to be the classical $q$%
-Bernstein-Kantorovich operators.\newline
\newline
First, we prove the following basic lemmas:\newline

\parindent=0mm\textbf{Lemma 2.1}. For $x\in \lbrack 0,1],~0<q<p\leq 1$

\begin{itemize}
\item[(i)] $K_{n}^{(p,q)}(1;x)=1$,

\item[(ii)] $K_{n}^{(p,q)}(t;x)=x+ \frac{p^n}{[2]_{p,q}[n]_{p,q}}$,

\item[(iii)] $K_{n}^{(p,q)}(t^2;x)=\frac{q}{p}\frac{[n-1]_{p,q}}{[n]_{p,q}}%
x^2+\Big(\frac{p^n(2q+p)}{[3]_{p,q}[n]_{p,q}}+\frac{p^{n-1}}{[n]_{p,q}}\Big)%
x+ \frac{p^{2n}}{[3]_{p,q}[n]_{p,q}^2}$,

\item[(iv)] $K_{n}^{(p,q)}\big{(}(t-x)^2;x\big{)}=\Big(\frac{q}{p}\frac{%
[n-1]_{p,q}}{[n]_{p,q}}-1\Big)x^2  +\Big(\frac{p^n(2q+p)}{[3]_{p,q}[n]_{p,q}}%
+\frac{p^{n-1}}{[n]_{p,q}}-\frac{2p^n}{[2]_{p,q}[n]_{p,q}}\Big)x+ \frac{%
p^{2n}}{[3]_{p,q}[n]_{p,q}^2}$.
\end{itemize}

\parindent=0mm\textbf{Proof}. (i)
\begin{equation*}
K_{n}^{(p,q)}(1;x)=\frac{[n]_{p,q}}{p^{\frac{n(n-1)}{2}}}\sum%
\limits_{k=0}^{n} \frac{b_{n,k}^{(p,q)}(x)}{p^{n-k} q^k} \int_{\frac{%
[k]_{p,q}}{p^{k-n-1}~[n]_{p,q}}}^{\frac{[k+1]_{p,q}}{p^{k-n}~[n]_{p,q}}%
}d_{p,q}t=1.
\end{equation*}

(ii)
\begin{eqnarray*}
K_{n}^{(p,q)}(t;x)&=&\frac{[n]_{p,q}}{p^{\frac{n(n-1)}{2}}}%
\sum\limits_{k=0}^{n} \frac{b_{n,k}^{(p,q)}(x)}{p^{n-k} q^k} \int_{\frac{%
[k]_{p,q}}{p^{k-n-1}~[n]_{p,q}}}^{\frac{[k+1]_{p,q}}{p^{k-n}~[n]_{p,q}}%
}t~d_{p,q}t \\
&=&\frac{1}{[2]_{p,q}[n]_{p,q}}\frac{1}{p^{\frac{n(n-1)}{2}}}%
\sum\limits_{k=0}^{n} \frac{b_{n,k}^{(p,q)}(x)}{p^{n-k} q^k}\bigg(\frac{%
[k+1]^2_{p,q}-p^2[k]^2_{p,q}}{p^{2k-2n}}\bigg).
\end{eqnarray*}%
Using $[k+1]_{p,q}={q^{k}+p[k]_{p,q}}$, we have
\begin{eqnarray*}
K_{n}^{(p,q)}(t;x)&=&\frac{1}{[2]_{p,q}[n]_{p,q}}\frac{1}{p^{\frac{n(n-1)}{2}%
}}\sum\limits_{k=0}^{n} \frac{b_{n,k}^{(p,q)}(x)(p^k+[2]_{p,q}[k]_{p,q})}{%
p^{n-k} q^k} \\
&=&\frac{1}{[2]_{p,q}[n]_{p,q}}\frac{p^n}{p^{\frac{n(n-1)}{2}}} \bigg(%
\sum\limits_{k=0}^{n}b_{n,k}^{(p,q)}(x)+[2]_{p,q}\sum%
\limits_{k=0}^{n}b_{n,k}^{(p,q)}(x)\frac{[k]_{p,q}}{p^k}\bigg) \\
&=&\frac{1}{[2]_{p,q}[n]_{p,q}}\frac{p^n}{p^{\frac{n(n-1)}{2}}} \bigg(p^{%
\frac{n(n-1)}{2}}+[2]_{p,q}\sum\limits_{k=0}^{n}\left[
\begin{array}{c}
n \\
k%
\end{array}%
\right]_{p,q}(x)_{p,q}^k(1-x)_{p,q}^{n-k}\frac{[k]_{p,q}}{p^k}\bigg) \\
&=&\frac{p^n}{[2]_{p,q}[n]_{p,q}}+\frac{1}{p^{\frac{n(n-3)}{2}}}%
\sum\limits_{k=0}^{n-1}\left[
\begin{array}{c}
n-1 \\
k%
\end{array}%
\right]_{p,q}(x)_{p,q}^{k+1}(1-x)_{p,q}^{n-k-1}\frac{1}{p^{k+1}} \\
&=&\frac{p^n}{[2]_{p,q}[n]_{p,q}}+\frac{1}{p^{\frac{(n-1)(n-2)}{2}}}%
\sum\limits_{k=0}^{n-1}\left[
\begin{array}{c}
n-1 \\
k%
\end{array}%
\right]_{p,q}(p^kx)(x)_{p,q}^{k}(1-x)_{p,q}^{n-k-1}\frac{1}{p^{k}} \\
&=&\frac{p^n}{[2]_{p,q}[n]_{p,q}}+\frac{x}{p^{\frac{(n-1)(n-2)}{2}}}%
\sum\limits_{k=0}^{n-1}\left[
\begin{array}{c}
n-1 \\
k%
\end{array}%
\right]_{p,q}(x)_{p,q}^{k}(1-x)_{p,q}^{n-k-1} \\
&=&\frac{p^n}{[2]_{p,q}[n]_{p,q}}+x.
\end{eqnarray*}%
(iii)
\begin{eqnarray*}
K_{n}^{(p,q)}(t^2;x)&=&\frac{[n]_{p,q}}{p^{\frac{n(n-1)}{2}}}%
\sum\limits_{k=0}^{n} \frac{b_{n,k}^{(p,q)}(x)}{p^{n-k} q^k} \int_{\frac{%
[k]_{p,q}}{p^{k-n-1}~[n]_{p,q}}}^{\frac{[k+1]_{p,q}}{p^{k-n}~[n]_{p,q}}%
}t^2~d_{p,q}t \\
&=&\frac{1}{[3]_{p,q}[n]^2_{p,q}}\frac{1}{p^{\frac{n(n-1)}{2}}}%
\sum\limits_{k=0}^{n} \frac{b_{n,k}^{(p,q)}(x)}{p^{n-k} q^k}\bigg(\frac{%
[k+1]^3_{p,q}-p^3[k]^3_{p,q}}{p^{3k-3n}}\bigg) \\
&=&\frac{p^{2n}}{[3]_{p,q}[n]^2_{p,q}}\frac{1}{p^{\frac{n(n-1)}{2}}}%
\sum\limits_{k=0}^{n}b_{n,k}^{(p,q)}(x) \bigg(1+(2q+p)\frac{[k]_{p,q}}{p^k}+%
\frac{[3]_{p,q}[k]^2_{p,q}}{p^{2k}}\bigg) \\
&=&\frac{p^{2n}}{[3]_{p,q}[n]^2_{p,q}}\bigg(1+\frac{(2q+p)}{p^{\frac{n(n-1)}{%
2}}}\sum\limits_{k=0}^{n}b_{n,k}^{(p,q)}(x)\frac{[k]_{p,q}}{p^k} +\frac{%
[3]_{p,q}}{p^{\frac{n(n-1)}{2}}}\sum\limits_{k=0}^{n}b_{n,k}^{(p,q)}(x)\frac{%
[k]^2_{p,q}}{p^{2k}}\bigg).
\end{eqnarray*}%
With the help of the previous calculations, we have
\begin{eqnarray*}
\frac{1}{p^{\frac{n(n-1)}{2}}}\sum\limits_{k=0}^{n}b_{n,k}^{(p,q)}(x)\frac{%
[k]_{p,q}}{p^k}&=&\frac{[n]_{p,q}}{p^n}x.
\end{eqnarray*}
And
\begin{eqnarray*}
\frac{1}{p^{\frac{n(n-1)}{2}}}\sum\limits_{k=0}^{n}b_{n,k}^{(p,q)}(x)\frac{%
[k]^2_{p,q}}{p^{2k}} &=& \frac{1}{p^{\frac{n(n-1)}{2}}}\sum\limits_{k=0}^{n}%
\left[
\begin{array}{c}
n \\
k%
\end{array}%
\right]_{p,q}(x)_{p,q}^k(1-x)_{p,q}^{n-k}\frac{[k]^2_{p,q}}{p^{2k}} \\
&=& \frac{[n]_{p,q}}{p^{\frac{n(n-1)}{2}}}\sum\limits_{k=0}^{n-1}\left[
\begin{array}{c}
n-1 \\
k%
\end{array}%
\right]_{p,q}(x)_{p,q}^{k+1}(1-x)_{p,q}^{n-k-1}\frac{[k+1]_{p,q}}{p^{2k+2}}.
\end{eqnarray*}
Using $[k+1]_{p,q}={p^{k}+q[k]_{p,q}}$, we have
\begin{align*}
& \frac{1}{p^{\frac{n(n-1)}{2}}}\sum\limits_{k=0}^{n}b_{n,k}^{(p,q)}(x)\frac{%
[k]^2_{p,q}}{p^{2k}} \\
&= \frac{1}{p^{\frac{n(n-1)}{2}}}\sum\limits_{k=0}^{n-1}\left[
\begin{array}{c}
n-1 \\
k%
\end{array}%
\right]_{p,q}(x)_{p,q}^{k+1}(1-x)_{p,q}^{n-k-1}\bigg(\frac{1}{p^{k+2}}+\frac{%
q[k]_{p,q}}{p^{2k+2}}\bigg) \\
& = \frac{1}{p^{\frac{n(n-1)}{2}}}\bigg\{\frac{p^{\frac{(n-1)(n-2)}{2}}x}{p^2%
}+q[n-1]_{p,q}\sum\limits_{k=0}^{n-2}\left[
\begin{array}{c}
n-2 \\
k%
\end{array}%
\right]_{p,q}(x)_{p,q}^{k+2}(1-x)_{p,q}^{n-k-2}\frac{1}{p^{2k+4}}\bigg\} \\
& = \frac{1}{p^{\frac{n(n-1)}{2}}}\bigg\{\frac{p^{\frac{(n-1)(n-2)}{2}}x}{p^2%
}+q[n-1]_{p,q}\sum\limits_{k=0}^{n-2}\left[
\begin{array}{c}
n-2 \\
k%
\end{array}%
\right]_{p,q}(p^{k+1}x)(p^{k}x)(x)_{p,q}^{k}(1-x)_{p,q}^{n-k-2}\frac{1}{%
p^{2k+4}}\bigg\} \\
&=\frac{1}{p^{\frac{n(n-1)}{2}}}\bigg\{\frac{p^{\frac{(n-1)(n-2)}{2}}x}{p^2}%
+q[n-1]_{p,q}\frac{p^{\frac{(n-2)(n-3)}{2}}x^2}{p^4}\bigg\} \\
&=\frac{1}{p^{n+1}}x+\frac{q[n-1]_{p,q}}{p^{2n+1}}x^2.
\end{align*}

Using the above equalities, we have
\begin{equation*}
K_{n}^{(p,q)}(t^2;x)=\frac{q}{p}\frac{[n-1]_{p,q}}{[n]_{p,q}}x^2+\bigg(\frac{%
p^n(2q+p)}{[3]_{p,q}[n]_{p,q}}+\frac{p^{n-1}}{[n]_{p,q}}\bigg)x+ \frac{p^{2n}%
}{[3]_{p,q}[n]_{p,q}^2}.
\end{equation*}

(iv) Using the linearity of the operators $K_{n}^{(p,q)}$, we have
\begin{align*}
& K_{n}^{(p,q)}\big{(}(t-x)^2;x\big{)} \\
& =K_{n}^{(p,q)}(t^2;x)-2xK_{n}^{(p,q)}(t;x)+x^2K_{n}^{(p,q)}(1;x) \\
&=\frac{q}{p}\frac{[n-1]_{p,q}}{[n]_{p,q}}x^2+\bigg(\frac{p^n(2q+p)}{%
[3]_{p,q}[n]_{p,q}}+\frac{p^{n-1}}{[n]_{p,q}}\bigg)x+ \frac{p^{2n}}{%
[3]_{p,q}[n]_{p,q}^2}-2x\bigg(x+\frac{p^n}{[2]_{p,q}[n]_{p,q}}\bigg)+x^2 \\
&=\bigg(\frac{q}{p}\frac{[n-1]_{p,q}}{[n]_{p,q}}-1\bigg)x^2+\bigg(\frac{%
p^n(2q+p)}{[3]_{p,q}[n]_{p,q}} +\frac{p^{n-1}}{[n]_{p,q}}-\frac{2p^n}{%
[2]_{p,q}[n]_{p,q}}\bigg)x+ \frac{p^{2n}}{[3]_{p,q}[n]_{p,q}^2}. \\
\end{align*}

\section{Main Results}

Let $C[a,b]$ be the linear space of all real valued continuous functions $f$
on $[a,b]$ and let $T$ be a linear operator which maps $C[a,b]$ into itself.
We say that $T$ is $positive$ if for every non-negative $f\in $ $C[a,b],$ we
have $T(f,x)\geq 0$ for all $x\in $ $[a,b]$ .

\parindent=8mm The classical Korovkin approximation theorem \cite{alt, pp,
sr} states as follows:

\parindent=8mm Let $(T_{n})$ be a sequence of positive linear operators from
$\mathcal{C}[a,b]$ into $C[a,b].$ Then $\lim_{n}\Vert T_{n}(f,x)-f(x)\Vert
_{C[a,b]}=0$, for all $f\in C[a,b]$ if and only if $\lim_{n}\Vert
T_{n}(f_{i},x)-f_{i}(x)\Vert _{C[a,b]}=0$, for $i=0,1,2$, where $f_{0}(x)=1,$
$f_{1}(x)=x$ and $f_{2}(x)=x^{2}.$\newline

\parindent=0mm\textbf{Theorem 3.1.} Let $0<q_{n}<p_{n}\leq 1$ such that $%
\lim\limits_{n\rightarrow \infty }p_{n}=1$ and $\lim\limits_{n\rightarrow
\infty }q_{n}=1$. Then for each $f\in C[0,1],~K_{n}^{(p_{n},q_{n})}(f;x)$
converges uniformly to $f$ on $[0,1]$.\newline
\newline
\parindent=0mm\textbf{Proof}. By the Korovkin Theorem it is sufficient to
show that
\begin{equation*}
\lim\limits_{n\rightarrow \infty }\Vert
K_{n}^{(p_{n},q_{n})}(t^{m};x)-x^{m}\Vert _{C[0,1]}=0,~~~m=0,1,2.
\end{equation*}%
By Lemma 2.1 (i), it is clear that
\begin{equation*}
\lim\limits_{n\rightarrow \infty }\Vert K_{n}^{(p_{n},q_{n})}(1;x)-1\Vert
_{C[0,1]}=0.
\end{equation*}%
Now, by Lemma 2.1 (ii)
\begin{eqnarray*}
|K_{n}^{(p_{n},q_{n})}(t;x)-x|&=& \frac{p_n^n}{[2]_{p_n,q_n}[n]_{p_n,q_n}}
\end{eqnarray*}
which yields
\begin{equation*}
\lim\limits_{n\rightarrow \infty }\Vert K_{n}^{(p_{n},q_{n})}(t;x)-x\Vert
_{C[0,1]}=0.
\end{equation*}%
Similarly,
\begin{align*}
& |K_{n}^{(p_{n},q_{n})}(t^2;x)-x^2| \\
& =\bigg{|}\bigg(\frac{q_n}{p_n}\frac{[n-1]_{p_n,q_n}}{[n]_{p_n,q_n}}-1\bigg)%
x^2+\bigg(\frac{p_n^n(2q_n+p_n)}{[3]_{p_n,q_n}[n]_{p_n,q_n}} +\frac{p_n^{n-1}%
}{[n]_{p_n,q_n}}\bigg)x+ \frac{p_n^{2n}}{[3]_{p_n,q_n}[n]_{p_n,q_n}^2}%
\bigg{|} \\
&\le \bigg(\frac{q_n}{p_n}\frac{[n-1]_{p_n,q_n}}{[n]_{p_n,q_n}}-1\bigg)x^2+%
\bigg(\frac{p_n^n(2q_n+p_n)}{[3]_{p_n,q_n}[n]_{p_n,q_n}} +\frac{p_n^{n-1}}{%
[n]_{p_n,q_n}}\bigg)x+ \frac{p_n^{2n}}{[3]_{p_n,q_n}[n]_{p_n,q_n}^2}.
\end{align*}
Taking maximum of both sides of the above inequality, we get
\begin{equation*}
\|K_{n}^{(p_{n},q_{n})}(t^2;x)-x^2\|\le\frac{q_n}{p_n}\frac{[n-1]_{p_n,q_n}}{%
[n]_{p_n,q_n}}-1+\frac{p_n^n(2q_n+p_n)}{[3]_{p_n,q_n}[n]_{p_n,q_n}} +\frac{%
p_n^{n-1}}{[n]_{p_n,q_n}}+ \frac{p_n^{2n}}{[3]_{p_n,q_n}[n]_{p_n,q_n}^2}
\end{equation*}
which concludes
\begin{equation*}
\lim\limits_{n\rightarrow \infty }\Vert
K_{n}^{(p_{n},q_{n})}(t^2;x)-x^2\Vert _{C[0,1]}=0.
\end{equation*}%
Thus the proof is completed.\newline

Now we will compute the rate of convergence in terms of modulus of
continuity.

\parindent=8mm Let $f\in C[0,1]$. The modulus of continuity of $f$ denoted
by $\omega(f,\delta )$ gives the maximum oscillation of $f$ in any interval
of length not exceeding $\delta >0$ and it is given by the relation
\begin{equation*}
\omega(f,\delta )=\sup\limits_{|x-y|\leq \delta }|f(x)-f(y)|,~~x,y\in
\lbrack 0,b].
\end{equation*}%
It is known that $\lim\limits_{\delta \rightarrow 0^+}\omega(f,\delta )=0$
for $f\in C[0,b]$ and for any $\delta >0$ one has
\begin{equation*}
|f(y)-f(x)|\leq \omega(f,\delta)\biggl{(}\frac{(y-x)^{2}}{\delta ^{2}}+1%
\biggl{)}.\eqno(3.1)
\end{equation*}

\parindent=0mm \textbf{Theorem 3.2}. If $f\in C[0,1]$, then
\begin{equation*}
\bigl{|}K_{n}^{(p,q)}(f;x)-f(x)\bigl{|}\leq 2\omega\big(f,\delta _{n}(x)\big)
\end{equation*}%
takes place, where $\delta _{n}(x)=\sqrt{K_{n}^{(p,q)}\big{(}(t-x)^2}$.
\newline

\parindent=0mm \textbf{Proof}. Since $K_{n}^{(p,q)}(1;x)=1$, we have
\begin{eqnarray*}
\bigl{|}K_{n}^{(p,q)}(f;x)-f(x)\bigl{|}&\leq&K_{n}^{(p,q)}\big{(}%
|f(t)-f(x)|;x\big{)} \\
&\leq& \frac{[n]_{p,q}}{p^{\frac{n(n-1)}{2}}}\sum\limits_{k=0}^{n} \frac{%
b_{n,k}^{(p,q)}(x)}{p^{n-k} q^k} \int_{\frac{[k]_{p,q}}{p^{k-n-1}~[n]_{p,q}}%
}^{\frac{[k+1]_{p,q}}{p^{k-n}~[n]_{p,q}}}|f(t)-f(x)|d_{p,q}t.
\end{eqnarray*}
In view of (3.1), we get
\begin{eqnarray*}
\bigl{|}K_{n}^{(p,q)}(f;x)-f(x)\bigl{|}&\leq&\bigg{\{}\frac{[n]_{p,q}}{p^{%
\frac{n(n-1)}{2}}}\sum\limits_{k=0}^{n} \frac{b_{n,k}^{(p,q)}(x)}{p^{n-k} q^k%
} \int_{\frac{[k]_{p,q}}{p^{k-n-1}~[n]_{p,q}}}^{\frac{[k+1]_{p,q}}{%
p^{k-n}~[n]_{p,q}}}\bigg{(}\frac{|t-x|^2}{\delta^2}+1\bigg{)}d_{p,q}t%
\bigg{\}}\omega(f,\delta) \\
&=&\Big{\{}\frac{1}{\delta^2}K_{n}^{(p,q)}\big{(}(t-x)^2;x\big{)}+1\Big{\}}%
\omega(f,\delta). \\
\end{eqnarray*}
Choosing $\delta=\delta _{n}(x)=\sqrt{K_{n}^{(p,q)}\big{(}(t-x)^2}$, we have
\begin{equation*}
\bigl{|}K_{n}^{(p,q)}(f;x)-f(x)\bigl{|}\leq 2\omega\big(f,\delta _{n}(x)\big)%
.
\end{equation*}
This completes the proof of the theorem.\newline

\parindent=8mmNow we give the rate of convergence of the operators $%
K_{n}^{(p,q)}$ in terms of the elements of the usual Lipschitz class $\text{%
Lip}_M(\alpha)$.\newline
$~~~~~~~~~~$Let $f\in C[0,1]$, $M>0$ and $0<\alpha\leq1$.We recall that $f$
belongs to the class $\text{Lip}_M(\alpha)$ if the inequality
\begin{equation*}
|f(t)-f(x)|\leq M|t-x|^\alpha,~~~(t,x\in[0,1])
\end{equation*}
is satisfied.\newline

%\newpage
\parindent=0mm\textbf{Theorem 3.3}. Let $0<q<p\leq1$. Then for each $f\in
\text{Lip}_M(\alpha)$ we have
\begin{equation*}
|K_{n}^{(p,q)}(f;x)-f(x)|\leq M\delta_n^{\alpha}(x),
\end{equation*}
where $\delta _{n}(x)=\sqrt{K_{n}^{(p,q)}\big{(}(t-x)^2;x\big{)}}.$\newline

\parindent=0mm\textbf{Proof}. By the monotonicity of the operators $%
K_{n}^{(p,q)}$, we can write
\begin{eqnarray*}
\bigl{|}K_{n}^{(p,q)}(f;x)-f(x)\bigl{|}&\leq&K_{n}^{(p,q)}\big{(}%
|f(t)-f(x)|;x\big{)} \\
&\leq& \frac{[n]_{p,q}}{p^{\frac{n(n-1)}{2}}}\sum\limits_{k=0}^{n} \frac{%
b_{n,k}^{(p,q)}(x)}{p^{n-k} q^k} \int_{\frac{[k]_{p,q}}{p^{k-n-1}~[n]_{p,q}}%
}^{\frac{[k+1]_{p,q}}{p^{k-n}~[n]_{p,q}}}|f(t)-f(x)|d_{p,q}t \\
&\leq& M\frac{[n]_{p,q}}{p^{\frac{n(n-1)}{2}}}\sum\limits_{k=0}^{n} \frac{%
b_{n,k}^{(p,q)}(x)}{p^{n-k} q^k} \int_{\frac{[k]_{p,q}}{p^{k-n-1}~[n]_{p,q}}%
}^{\frac{[k+1]_{p,q}}{p^{k-n}~[n]_{p,q}}}|t-x|^\alpha~d_{p,q}t.
\end{eqnarray*}

Now applying the H\"{o}lder's inequality for the sum with $p_1=\frac2{\alpha}
$ and $p_2=\frac2{2-\alpha}$ and taking into consideration Lemma 2.1(i) and
Lemma 2.2(ii), we have
\begin{eqnarray*}
|K_{n}^{(p,q)}(f;x)-f(x)| &\leq& M\sum\limits_{k=0}^{n}\Bigg{\{}\frac{%
[n]_{p,q}}{p^{\frac{n(n-1)}{2}}} \frac{b_{n,k}^{(p,q)}(x)}{p^{n-k} q^k}
\int_{\frac{[k]_{p,q}}{p^{k-n-1}~[n]_{p,q}}}^{\frac{[k+1]_{p,q}}{%
p^{k-n}~[n]_{p,q}}}(t-x)^2 d_{p,q}t \Bigg{\}}^{\frac\alpha2} \\
&&\times \Bigg{\{}\frac{[n]_{p,q}}{p^{\frac{n(n-1)}{2}}} \frac{%
b_{n,k}^{(p,q)}(x)}{p^{n-k} q^k} \int_{\frac{[k]_{p,q}}{p^{k-n-1}~[n]_{p,q}}%
}^{\frac{[k+1]_{p,q}}{p^{k-n}~[n]_{p,q}}}1 d_{p,q}t \Bigg{\}}^{\frac{2-\alpha%
}2} \\
&\leq& M\Bigg{\{}\frac{[n]_{p,q}}{p^{\frac{n(n-1)}{2}}} \sum\limits_{k=0}^{n}%
\frac{b_{n,k}^{(p,q)}(x)}{p^{n-k} q^k} \int_{\frac{[k]_{p,q}}{%
p^{k-n-1}~[n]_{p,q}}}^{\frac{[k+1]_{p,q}}{p^{k-n}~[n]_{p,q}}}(t-x)^2
d_{p,q}t \Bigg{\}}^{\frac\alpha2} \\
&&\times \Bigg{\{}\frac{[n]_{p,q}}{p^{\frac{n(n-1)}{2}}} \sum%
\limits_{k=0}^{n}\frac{b_{n,k}^{(p,q)}(x)}{p^{n-k} q^k} \int_{\frac{[k]_{p,q}%
}{p^{k-n-1}~[n]_{p,q}}}^{\frac{[k+1]_{p,q}}{p^{k-n}~[n]_{p,q}}}1 d_{p,q}t %
\Bigg{\}}^{\frac{2-\alpha}2} \\
&=&M\bigl{\{}K_{n}^{(p,q)}\bigl{(}(t-x)^2;x\bigl{)}\bigl{\}}^\frac{\alpha}{2}%
.
\end{eqnarray*}
Choosing $\delta^2(x)=\delta_n^2(x)=K_{n}^{(p,q)}\big{(}(t-x)^2;x\big{)}$,
we arrive at our desired result. %\newline
\newline

\parindent=8mm Next, we prove the local approximation property for the
operators $K_{n}^{(p,q)}$. The Peetre's $K$-functional is defined by \newline
\begin{equation*}
K_{2}(f,\delta )=\inf [\{\Vert f-g\Vert +\delta \Vert g^{\prime \prime
}\Vert \}:g\in W^{2}],
\end{equation*}%
where%
\begin{equation*}
W^{2}=\{g\in C[0,1]:g^{\prime },g^{\prime \prime }\in C[0,1]\}.
\end{equation*}

By \cite{dl}, there exists a positive constant $C>0$ such that $%
K_{2}(f,\delta )\leq C\omega_{2}(f,\sqrt{\delta}),~\delta >0$, where the
second order modulus of continuity is given by
\begin{equation*}
\omega_{2}(f,\sqrt{\delta})=\sup\limits_{0<h\leq \sqrt{\delta}%
}\sup\limits_{x\in \lbrack 0,1]}\mid f(x+2h)-2f(x+h)+f(x)\mid .
\end{equation*}%
%
%
%Also for $f\in [0,1]$ the usual modulus of continuity is given by
%\begin{equation*}
%\omega(f,\delta )=\sup\limits_{0<h\leq \delta ^{\frac{1}{2}%
%}}\sup\limits_{x\in \lbrack 0,1]}\mid f(x+h)-f(x)\mid .
%\end{equation*}%
%\newline

\parindent=0mm \textbf{Theorem 3.4}. Let $f\in C\lbrack 0,1]$ and $0<q<p\leq
1 $. Then for all $n\in \mathbb{N}$, there exists an absolute constant $C>0$
such that
\begin{equation*}
\big|K_n^{(p,q)}(f;x)-f(x)\big|\leq C\omega_2\big(f,\delta_n(x)\big)+\omega%
\big(f,\alpha_n(x)\big)
\end{equation*}
where%
\begin{equation*}
\delta _{n}(x)=\sqrt{K_{n}^{(p,q)}\left( (t-x)^{2};x\right)+\frac{p^{2n}}{%
[2]_{p,q}^2[n]_{p,q}^2}}\mbox{~~and~~} \alpha_n=\frac{p^{n}}{%
[2_{p,q}][n]_{p,q}}.
\end{equation*}
\newline

\parindent=0mm \textbf{Proof}. For $x\in[0,1]$, we consider the auxiliary
operators $K_n^*$ defined by
\begin{equation*}
K_n^*(f;x)=K_{n}^{(p,q)}(f;x)+f(x)-f\bigg(x+\frac{p^{n}}{[2_{p,q}][n]_{p,q}}%
\bigg).
\end{equation*}
From Lemma 2.1, we observe that the operators $K_n^*(f;x)$ are linear and
reproduce the linear functions. Hence
\begin{eqnarray*}
K_n^*(1;x) &=& K_n^{(p,q)}(1;x)+1-1=1, \\
K_n^*(t;x)&=& K_n^{(p,q)}(t;x)+x-\bigg(x+\frac{p^{n}}{[2_{p,q}][n]_{p,q}}%
\bigg)=x, \\
\mbox{so~~~} K_n^*(t-x;x)&=& K_n^*(t;x)-xK_n^*(1;x)=0.
\end{eqnarray*}

Let $x\in[0,1]$ and $g\in C_B^2[0,1]$. Using the Taylor's formula
\begin{equation*}
g(t)=g(x)+g^{\prime }(x)(t-x)+\int_{x}^{t}(t-u)~g^{\prime \prime }(u)~du.
\end{equation*}%
\newline
Applying $K_n^*$ to both sides of the above equation, we have
\begin{eqnarray*}
K_n^*(g;x)-g(x) &=& K_n^*\big((t-x)g^\prime(x);x\big)+K_n^*\bigg(%
\int_x^t(t-u)g^{\prime\prime}(u)du;x\bigg) \\
&=& g^\prime(x)K_n^*\big((t-x);x\big)+K_n^{(p,q)}\bigg(\int_x^t(t-u)g^{%
\prime\prime}(u)du;x\bigg) \\
&&-\int_x^{x+\frac{p^{n}}{[2_{p,q}][n]_{p,q}}} \bigg(x+\frac{p^{n}}{%
[2_{p,q}][n]_{p,q}}-u\bigg)g^{\prime\prime}(u)du \\
&=&K_n^{(p,q)}\bigg(\int_x^t(t-u)g^{\prime\prime}(u)du;x\bigg) \\
&&-\int_x^{x+\frac{p^{n}}{[2_{p,q}][n]_{p,q}}} \bigg(x+\frac{p^{n}}{%
[2_{p,q}][n]_{p,q}}-u\bigg)g^{\prime\prime}(u)du.
\end{eqnarray*}
On the other hand, since
\begin{equation*}
\bigg|\int_x^t(t-u)g^{\prime\prime}(u)du\bigg|\leq\int_x^t|t-u||g^{\prime%
\prime}(u)|du\leq\|g^{\prime\prime}\|\int_x^t|t-u|du\leq(t-x)^2\|g^{\prime%
\prime}\|
\end{equation*}
and
\begin{eqnarray*}
\Bigg|\int_x^{x+\frac{p^{n}}{[2_{p,q}][n]_{p,q}}} \bigg(x+\frac{p^{n}}{%
[2_{p,q}][n]_{p,q}}-u\bigg)g^{\prime\prime}(u)du\Bigg| \leq \frac{p^{2n}}{%
[2]_{p,q}^2[n]_{p,q}^2}\|g^{\prime\prime}\|.
\end{eqnarray*}
We conclude that
\begin{eqnarray*}
\Big|K_n^*(g;x)-g(x)\Big| &=& \bigg|K_n^{(p,q)}\bigg(\int_x^t(t-u)g^{\prime%
\prime}(u)du;x\bigg) \\
&&-\int_x^{x+\frac{p^{n}}{[2_{p,q}][n]_{p,q}}} \bigg(x+\frac{p^{n}}{%
[2_{p,q}][n]_{p,q}}-u\bigg)g^{\prime\prime}(u)du\bigg| \\
&\leq& \Vert g^{\prime \prime }\Vert K_{n}^{(p,q)}\left( (t-x)^{2};x\right)+%
\frac{p^{2n}}{[2]_{p,q}^2[n]_{p,q}^2}\|g^{\prime\prime}\| \\
&=&\delta_n^2(x)\Vert g^{\prime \prime }\Vert.
\end{eqnarray*}
Now, taking into account boundedness of $K_n^*$, we have
\begin{equation*}
\big|K_n^*(f;x)\big|\leq\big|K_n^{(p,q)}(f;x)\big|+2\|f\|\leq3\|f\|.
\end{equation*}
Therefore
\begin{eqnarray*}
\big|K_n^{(p,q)}(f;x)-f(x)\big| &\leq& \big|K_n^*(f;x)-f(x)\big|+\bigg|f(x)-f%
\bigg(x+\frac{p^{n}}{[2]_{p,q}[n]_{p,q}}\bigg)\bigg| \\
&\leq& \big|K_n^*(f-g;x)-(f-g)(x)\big| \\
&&+\bigg|f(x)-f\bigg(x+\frac{p^{n}}{[2]_{p,q}[n]_{p,q}}\bigg)\bigg| +\big|%
K_n^*(g;x)-g(x)\big| \\
&\leq& \big|K_n^*(f-g;x)\big|+\big|(f-g)(x)\big| \\
&&+\bigg|f(x)-f\bigg(x+\frac{p^{n}}{[2]_{p,q}[n]_{p,q}}\bigg)\bigg| +\big|%
K_n^*(g;x)-g(x)\big| \\
&\leq&4\|f-g\|+\omega\bigg(f,\frac{p^{n}}{[2]_{p,q}[n]_{p,q}}\bigg) %
+\delta_n^2(x)\Vert g^{\prime \prime }\Vert.
\end{eqnarray*}
Hence, taking the infimum on the right-hand side over all $g\in W^2$, we
have the following result
\begin{equation*}
\big|K_n^{(p,q)}(f;x)-f(x)\big|\leq4K_2\big(f,\delta_n^2(x)\big)+\omega\big(%
f,\alpha_n(x)\big).
\end{equation*}
In view of the property of $K$-functional, we get
\begin{equation*}
\big|K_n^{(p,q)}(f;x)-f(x)\big|\leq C\omega_2\big(f,\delta_n(x)\big)+\omega%
\big(f,\alpha_n(x)\big).
\end{equation*}

\parindent=0mm This completes the proof of the theorem.

\section{Graphical Examples}

With the help of Matlab, we show comparisons and some illustrative graphics
\cite{ma1} for the convergence of operators to the function $f (x) =
1+sin(7x)$.\newline

From Fig. 1 it can be observed that as the value of $q$ and $p$ approaches
towards $1$ provided $0 < q < p \leq 1,$ $(p,q)$-Bernstein-Kantorovich
operators converge towards the function. The parameters $q$ and $p$ adds
flexibility in approximation of functions by positive linear operators.%
\newline

In comparison to Fig. 1 as the value the $n$ increases, operators given by
(2.1) converge towards the function which is shown in figure $(2) $.\newline

Similarly for different values of parameters $p,q,n$ convergence of
operators to the function is shown in Fig. 3 and Fig. 4.\newline

Thus in comparison to $q$-analogoue of Bernstein-Kantorovich operators, our
generalization gives us more flexibility for the convergence of operators to
a function.

\begin{figure}[tbp]
\begin{center}
%\centering
\begin{minipage}{.3\textwidth}
%\centering
\includegraphics[height=5cm, width=7cm]{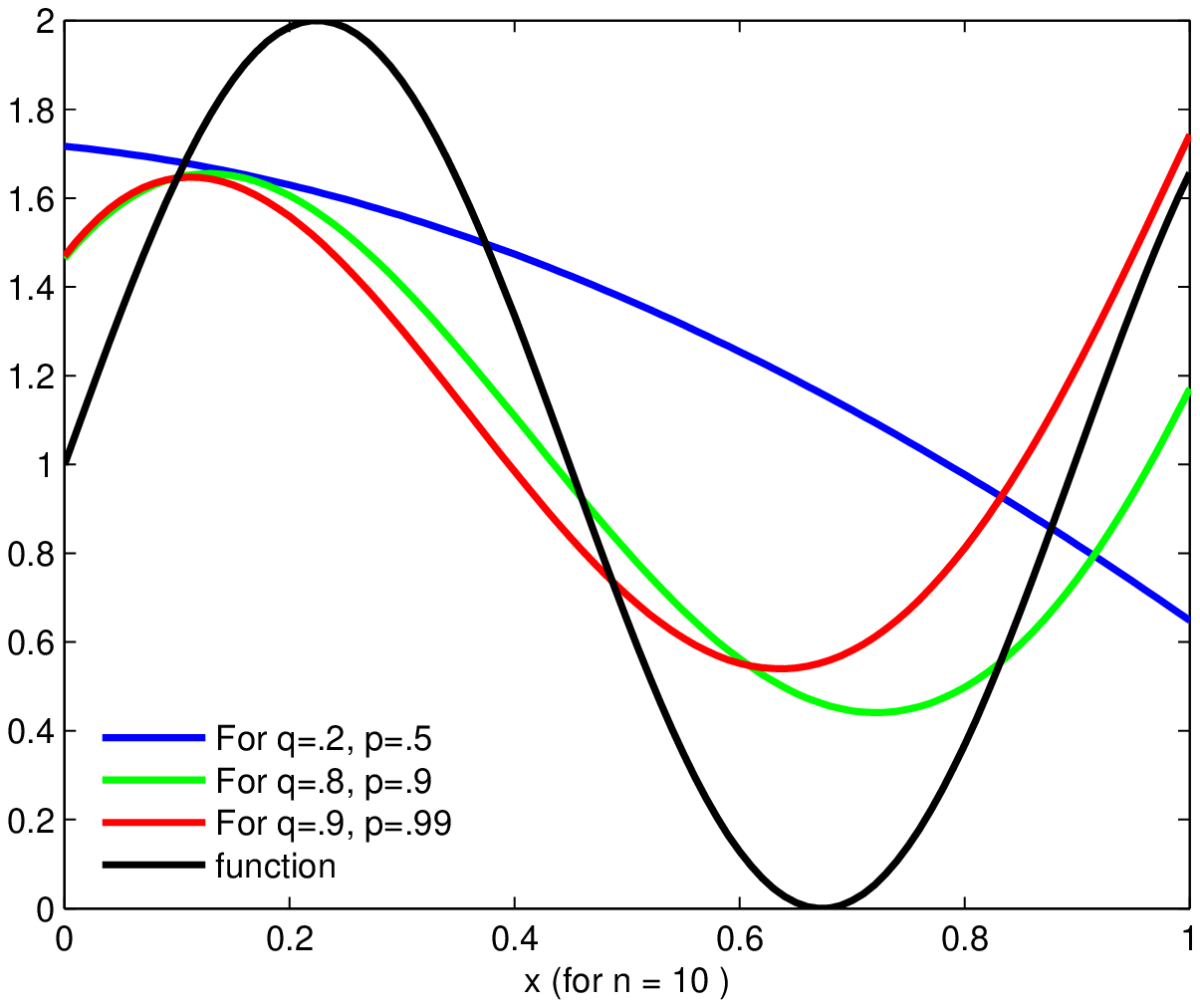}
\caption{}
\label{fig:test1}
\end{minipage}\hspace{3cm}
\begin{minipage}{.3\textwidth}
%\centering
\includegraphics[height=5cm, width=7cm]{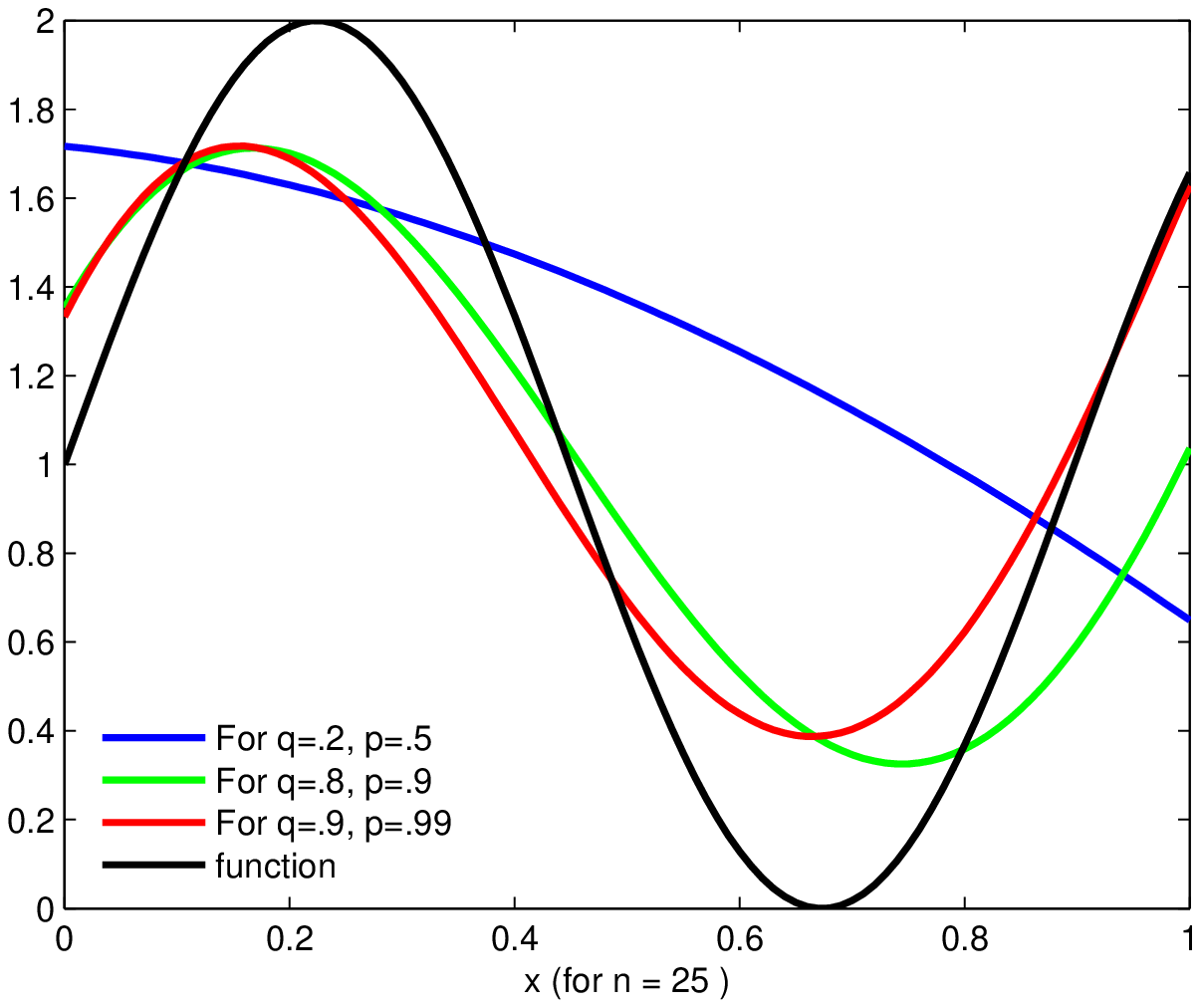}
\caption{}
\label{fig:test2}
\end{minipage}
\end{center}
\end{figure}

\begin{figure}[tbp]
\begin{center}
%\centering
\begin{minipage}{.3\textwidth}
%\centering
\includegraphics[height=5cm, width=7cm]{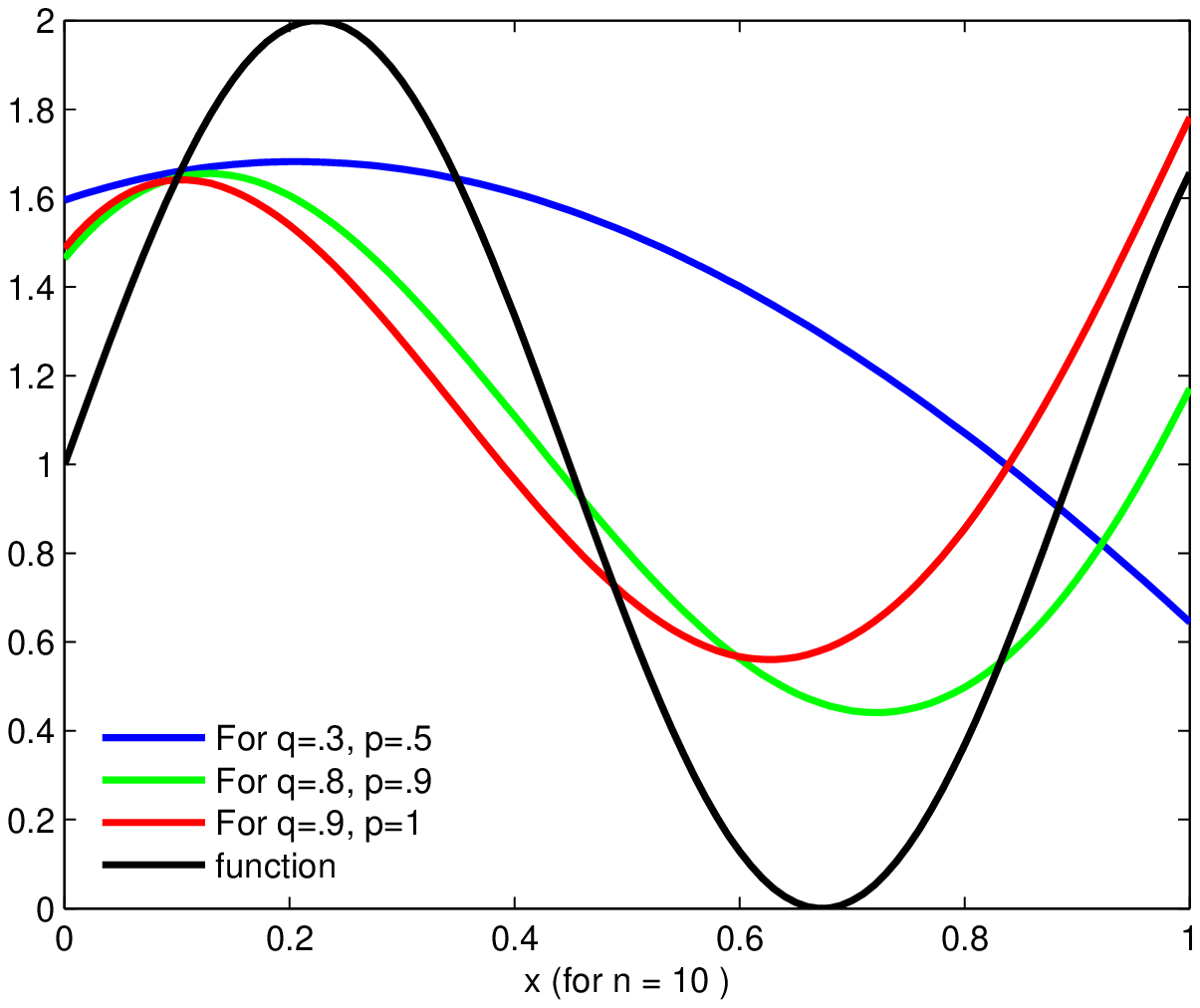}
\caption{}
\label{fig:test1}
\end{minipage}\hspace{3cm}
\begin{minipage}{.3\textwidth}
%\centering
\includegraphics[height=5cm, width=7cm]{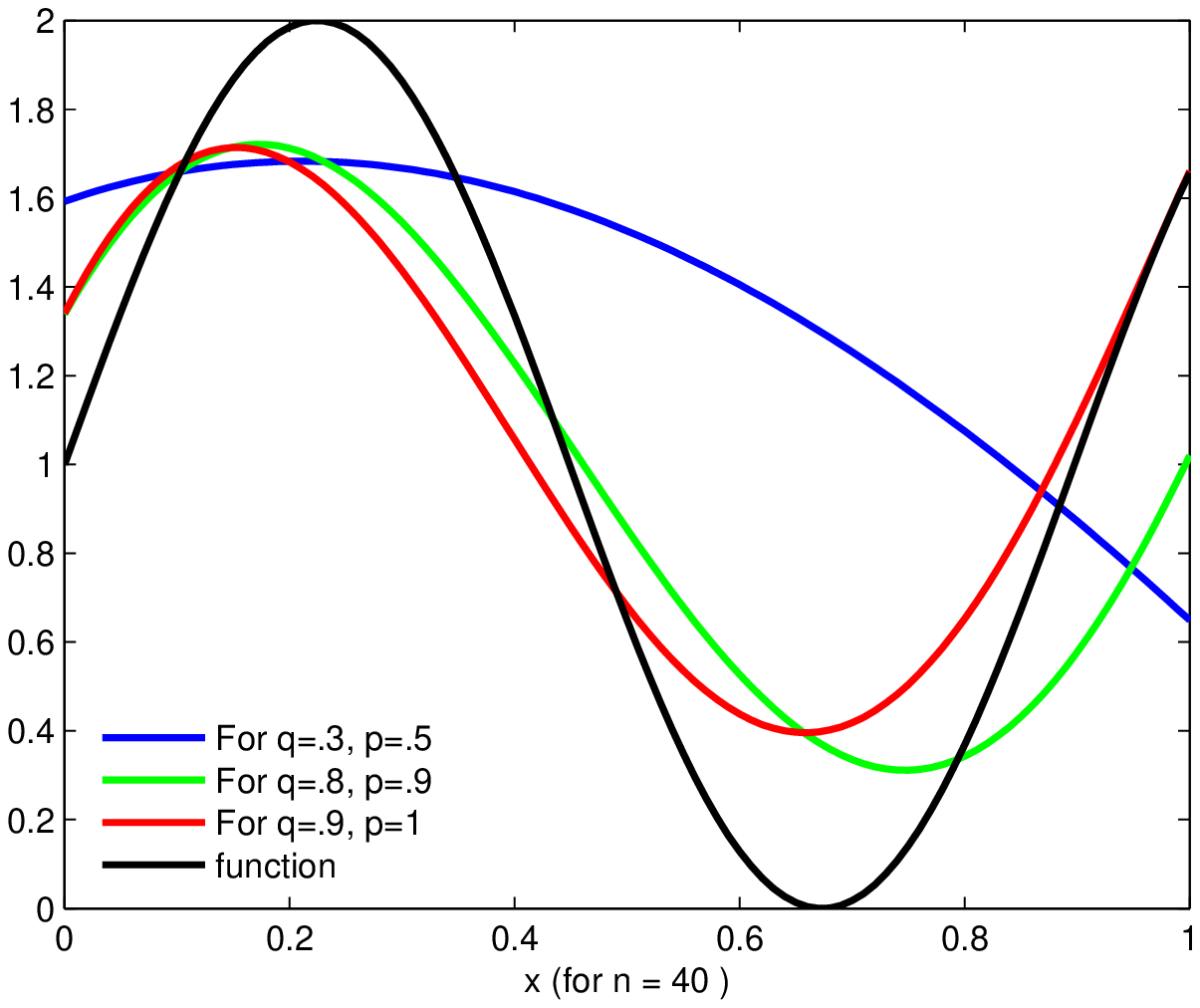}
\caption{}
\label{fig:test2}
\end{minipage}
\end{center}
\end{figure}

%\begin{figure*}[htb!]
%\begin{center}
%\includegraphics[height=5cm, width=7cm]{berst-kant5.eps}
%\end{center}
%\caption{}
%\end{figure*}
%
%
%\begin{figure*}[htb!]
%\begin{center}
%\includegraphics[height=5cm, width=7cm]{berst-kant3.eps}
%\end{center}
%\caption{}
%\end{figure*}
%
%
%
%
%\begin{figure*}[htb!]
%\begin{center}
%\includegraphics[height=5cm, width=7cm]{berst-kant6.eps}
%\end{center}
%\caption{}
%\end{figure*}
%\begin{figure*}[htb!]
%\begin{center}
%\includegraphics[height=5cm, width=7cm]{berst-kant7.eps}
%\end{center}
%\caption{}
%\end{figure*}

\end{document}